\documentclass[12pt]{article}%
\usepackage{amsmath}
\usepackage{amsfonts}
\usepackage{amssymb}
\usepackage{graphicx}
\usepackage{amssymb}%
\newtheorem{theorem}{Theorem}

\begin{document}
\title{Asymptotic analysis of the Askey-scheme II: from Charlier to Hermite}
\author{Diego Dominici \thanks{e-mail: dominicd@newpaltz.edu}\\Department of Mathematics\\State University of New York at New Paltz\\75 S. Manheim Blvd. Suite 9\\New Paltz, NY 12561-2443\\USA\\Phone: (845) 257-2607\\
Fax: (845) 257-3571
}

\maketitle

\begin{abstract}
We analyze the Hermite polynomials $H_{n}(\xi)$ and their zeros asymptotically
as $n\rightarrow\infty,$ using the limit relation between the Charlier and
Hermite polynomials. Our formulas involve some special functions and they
yield very accurate approximations.
\end{abstract}


MSC-class: 33C45 (Primary) 34E05, 33C10 (Secondary)


\section{Introduction}

The Hermite polynomials $H_{n}(x)$ are defined by \cite{MR1376370}
\begin{equation}
H_{n}(x)=n!%
{\displaystyle\sum\limits_{k=0}^{\left\lfloor \frac{n}{2}\right\rfloor }}
\frac{\left(  -1\right)  ^{k}}{k!(n-2k)!}\left(  2x\right)  ^{n-2k}
\label{Cndef}%
\end{equation}
for $n=0,1,\ldots.$ They satisfy the orthogonality condition \cite{atlas}
\[%
{\displaystyle\int\limits_{-\infty}^{\infty}}
e^{-x^{2}}H_{m}(x)H_{n}(x)dx=\sqrt{\pi}2^{n}n!\delta_{nm}%
\]
and the reflection formula
\begin{equation}
H_{n}(-x)=(-1)^{n}H_{n}(x). \label{reflec}%
\end{equation}
The Hermite polynomials are special cases of the parabolic cylinder function
$U(a,z),$
\[
H_{n}(x)=2^{\frac{n}{2}}\exp\left(  \frac{x^{2}}{2}\right)  U\left(
-n-\frac{1}{2},\sqrt{2}x\right) ,
\]
which was analyzed by Nico Temme in \cite{MR2076056}, \cite{MR1780050} and
\cite{MR1993339}.

The Hermite polynomials have been extensively studied since the pioneer
article of C. Hermite \cite{hermpoly} in 1864 (they were previously considered
by Fourier and Chebyshev). They have many applications in the physical
sciences and are particularly important in the quantum mechanical treatment of
the harmonic oscillator \cite{MR810399} (see also \cite{MR1972755},
\cite{MR1736478} and \cite{MR1343548} for some extensions). We refer the
interested reader to \cite{MR0058756} and \cite{MR0372517} for further
properties and references.

There are several families of orthogonal polynomials which have asymptotic
approximations in terms of $H_{n}(x).$ Some cases studied by Nico Temme
include the Gegenbauer \cite{MR1052447}, Laguerre \cite{MR2068115},
\cite{MR1036513}, Tricomi-Carlitz and Jacobi polynomials \cite{MR1803886}. He
also considered the asymptotic representations of other families of
polynomials such as the generalized Bernoulli, Euler, Bessel and Buchholz
polynomials in \cite{MR1723073}.

A rich source of asymptotic relations between the $H_{n}(x)$ and other
polynomials \cite{MR1985821}, \cite{MR1805994}, \cite{MR1858318} is provided
by the Askey-scheme of hypergeometric orthogonal polynomials
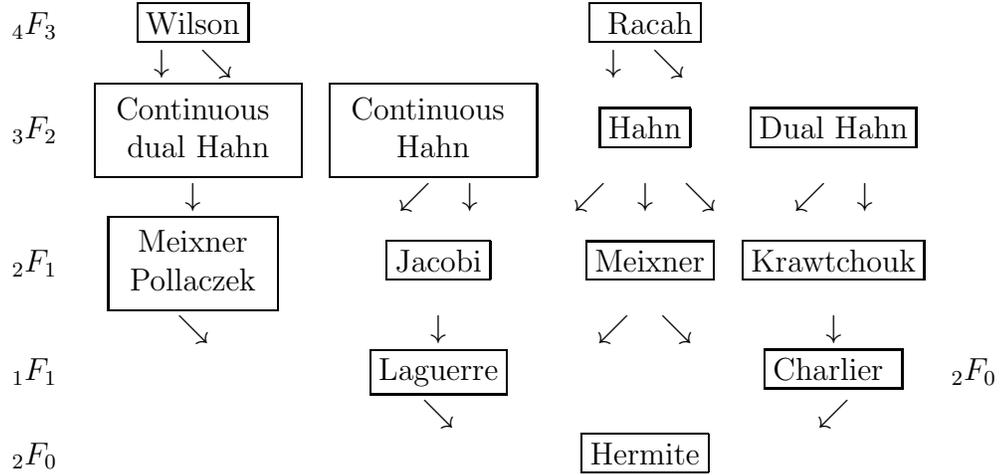
\cite{koekoek94askeyscheme}:

$%
\begin{array}
[c]{ccccccccc}%
_{4}F_{3} & \fbox{Wilson} & \  & \fbox{\ Racah} &  &  &  &  & \\
& \downarrow\quad\searrow &  & \downarrow\quad\searrow &  &  &  &  & \\
_{3}F_{2} & \ \fbox{$%
\begin{array}
[c]{c}%
\text{Continuous }\\
\text{dual Hahn}%
\end{array}
$} & \fbox{$%
\begin{array}
[c]{c}%
\text{Continuous }\\
\text{Hahn}%
\end{array}
$}\  & \fbox{Hahn} & \fbox{Dual Hahn} &  &  &  & \\
& \downarrow & \swarrow\quad\downarrow & \swarrow\quad\downarrow\quad\searrow
& \swarrow\quad\downarrow &  &  &  & \\
_{2}F_{1} & \fbox{$%
\begin{array}
[c]{c}%
\text{Meixner}\\
\text{Pollaczek}%
\end{array}
$} & \fbox{Jacobi} & \ \fbox{Meixner} & \fbox{Krawtchouk} &  &  &  & \\
& \searrow & \downarrow & \swarrow\quad\searrow & \downarrow &  &  &  & \\
_{1}F_{1} &  & \fbox{Laguerre} &  & \fbox{Charlier\ } & _{2}F_{0} &  &  & \\
&  & \searrow &  & \swarrow &  &  &  & \\
_{2}F_{0} &  &  & \fbox{Hermite} &  &  &  &  &
\end{array}
$

where the arrows indicate limit relations between the polynomials.

In particular, the limit relation between the Charlier polynomials
$C_{n}^{(a)}(x)$ and the Hermite polynomials is given by%
\begin{equation}
\underset{a\rightarrow\infty}{\lim}\left(  -1\right)  ^{n}\left(  2a\right)
^{\frac{n}{2}}C_{n}^{(a)}\left(  a+x\sqrt{2a}\right)  =H_{n}(x).\label{limit}%
\end{equation}

In this article we investigate the asymptotic behavior of $H_{n}(x)$ as
$n\rightarrow\infty,$ using (\ref{limit}) and the asymptotic results on the
Charlier polynomials derived in \cite{charlier}. We believe that our method
provides a useful approach to asymptotic analysis and could be used for other
families of polynomials of the Askey scheme.

\section{Previous results}

We define the Charlier polynomials by
\begin{equation}
\,C_{n}^{(a)}(x)=\,_{2}F_{0}\left(  \left.
\begin{array}
[c]{c}%
-n,-x\\
-
\end{array}
\right\vert -\frac{1}{a}\right)  ,\quad n=0,1,\ldots\label{Charlier}%
\end{equation}
with $a>0.$ The following results were derived in \cite{charlier}.

\begin{theorem}
As $n\rightarrow\infty,$ $C_{n}^{(a)}(x)$ admits the following asymptotic
approximations, with%
\begin{equation}
\Omega^{\pm}=\left(  \sqrt{n}\pm\sqrt{a}\right)  ^{2}. \label{Omega}%
\end{equation}

\begin{enumerate}
\item $n=O(1).$%
\begin{equation}
C_{n}^{(a)}\simeq\left(  1-\frac{x}{a}\right)  ^{n} \label{c1}%
\end{equation}

\item $x<\Omega^{-},\quad0<n<a.$
\begin{equation}
C_{n}^{(a)}\sim F_{3}(x)=\exp\left[  \Psi_{3}(x)\right]  L_{3}(x), \label{c2}%
\end{equation}
where
\begin{equation}
\Psi_{3}(x)=x\ln\left(  \frac{a+x-n+\Delta}{2a}\right)  +n\ln\left(
\frac{a-x+n+\Delta}{2a}\right)  +\frac{1}{2}\left(  a-x-n+\Delta\right)
\label{psi3}%
\end{equation}
and%
\begin{equation}
L_{3}(x)=\sqrt{\frac{a-x-n+\Delta}{2\Delta}}, \label{l3}%
\end{equation}
with%
\begin{equation}
\Delta(x)=\sqrt{a^{2}-2a(x+n)+(x-n)^{2}}. \label{delta}%
\end{equation}

\item $\Omega^{+}<x.$
\begin{equation}
C_{n}^{(a)}\sim F_{4}(x)=\left(  -1\right)  ^{n}\exp\left[  \Psi
_{4}(x)\right]  L_{4}(x), \label{c4}%
\end{equation}
where%
\begin{equation}
\Psi_{4}(x)=x\ln\left(  \frac{a+x-n-\Delta}{2a}\right)  +n\ln\left(
\frac{x-a-n+\Delta}{2a}\right)  +\frac{1}{2}\left(  a-x-n+\Delta\right)
\label{psi4}%
\end{equation}
and%
\begin{equation}
L_{4}(x)=\sqrt{\frac{x-a+n+\Delta}{2\Delta}}. \label{l4}%
\end{equation}

\item $x\approx\Omega^{-},\ 0<n<a.$%
\begin{align}
C_{n}^{(a)}  &  \sim\sqrt{2\pi}\left(  \frac{n}{a}\right)  ^{\frac{1}{6}%
}\left(  \sqrt{a}-\sqrt{n}\right)  ^{\frac{1}{3}}\mathrm{Ai}\left[  \left(
\frac{n}{a}\right)  ^{\frac{1}{6}}\frac{\left(  \Omega^{-}-x\right)  }{\left(
\sqrt{a}-\sqrt{n}\right)  ^{\frac{2}{3}}}\right] \label{c8}\\
&  \times\exp\left[  \frac{1}{2}n\ln\left(  \frac{n}{a}\right)  +x\ln\left(
1-\sqrt{\frac{n}{a}}\right)  +\sqrt{an}-\sqrt{n}\right]  ,\nonumber
\end{align}
where $\mathrm{Ai}\left(  \cdot\right)  $ is the Airy function.

\item $\Omega^{-}<x<\Omega^{+}.$%
\begin{equation}
C_{n}^{(a)}\sim F_{10}(x)=F_{3}(x)+F_{4}(x). \label{c10}%
\end{equation}

\item $x\approx\Omega^{+}.$%
\begin{align}
C_{n}^{(a)}  &  \sim\sqrt{2\pi}\left(  \frac{n}{a}\right)  ^{\frac{1}{6}%
}\left(  \sqrt{a}+\sqrt{n}\right)  ^{\frac{1}{3}}\mathrm{Ai}\left[  \left(
\frac{n}{a}\right)  ^{\frac{1}{6}}\frac{\left(  x-\Omega^{+}\right)  }{\left(
\sqrt{a}+\sqrt{n}\right)  ^{\frac{2}{3}}}\right] \label{c11}\\
&  \times\left(  -1\right)  ^{n}\exp\left[  \frac{1}{2}n\ln\left(  \frac{n}%
{a}\right)  +x\ln\left(  1+\sqrt{\frac{n}{a}}\right)  -\sqrt{an}-\sqrt
{n}\right]  .\nonumber
\end{align}

\end{enumerate}
\end{theorem}

\section{Limit analysis}

From (\ref{Omega}) we have, as $a\rightarrow\infty,$
\begin{equation}
\frac{\Omega^{\pm}-a}{\sqrt{2a}}\rightarrow\pm\sqrt{2n}. \label{limitOmega}%
\end{equation}
Hence, the six regions of Theorem 1 transform into the following regions:

\begin{enumerate}
\item Region I: $n=O(1).$

Setting
\begin{equation}
x=a+\xi\sqrt{2a} \label{eta}%
\end{equation}
in (\ref{c1}) and using (\ref{limit}), we get%
\begin{equation}
H_{n}(\xi)\simeq\left(  2\xi\right)  ^{n}. \label{h1}%
\end{equation}
The formula above is exact for $n=0,1$ and is a good approximation when
$\xi\gg n.$

\item Region II: $\xi<-\sqrt{2n}.$

From (\ref{delta}) and (\ref{eta}) we have
\begin{equation}
\Delta\sim\sqrt{2}\sigma\sqrt{a}-\frac{\xi n}{\sigma},\quad a\rightarrow
\infty, \label{delta1}%
\end{equation}
with%
\begin{equation}
\sigma=\sqrt{\xi^{2}-2n}. \label{sigma}%
\end{equation}
Using (\ref{delta1}) in (\ref{psi3}) and (\ref{l3}), we get%
\begin{equation}
\frac{n}{2}\ln(2a)+\Psi_{3}(x)\sim\Phi_{1}\left(  \xi\right)  \equiv
n\ln\left(  \sigma-\xi\right)  +\frac{\xi^{2}+\sigma\xi-n}{2}, \label{psi3lim}%
\end{equation}%
\begin{equation}
L_{3}(x)\sim U_{1}\left(  \xi\right)  =\sqrt{\frac{1}{2}\left(  1-\frac{\xi
}{\sigma}\right)  ,} \label{l3lim}%
\end{equation}
as $a\rightarrow\infty.$ Thus, from (\ref{limit}) we have%
\begin{equation}
H_{n}(\xi)\sim\Lambda_{1}\left(  \xi\right)  \equiv\left(  -1\right)  ^{n}%
\exp\left[  \Phi_{1}\left(  \xi\right)  \right]  U_{1}\left(  \xi\right)  .
\label{H2}%
\end{equation}

\item Region III: $\xi>\sqrt{2n}.$

Using (\ref{delta1}) in (\ref{psi4}) and (\ref{l4}), we obtain%
\begin{equation}
\frac{n}{2}\ln(2a)+\Psi_{4}(x)\sim\Phi_{2}\left(  \xi\right)  \equiv
n\ln\left(  \sigma+\xi\right)  +\frac{\xi^{2}-\sigma\xi-n}{2}, \label{psi4lim}%
\end{equation}%
\begin{equation}
L_{4}(x)\sim U_{2}\left(  \xi\right)  =\sqrt{\frac{1}{2}\left(  1+\frac{\xi
}{\sigma}\right)  ,} \label{l4lim}%
\end{equation}
as $a\rightarrow\infty.$ Hence,%
\begin{equation}
H_{n}(\xi)\sim\Lambda_{2}\left(  \xi\right)  \equiv\exp\left[  \Phi_{2}\left(
\xi\right)  \right]  U_{2}\left(  \xi\right)  . \label{H3}%
\end{equation}
Note that
\[
\Lambda_{1}\left(  -\xi\right)  =\left(  -1\right)  ^{n}\Lambda_{2}\left(
\xi\right)  ,
\]
as one would expect from (\ref{reflec}).

\item Region IV: $\xi\approx-\sqrt{2n}.$

Using (\ref{eta}) in (\ref{c8}) we have, as $a\rightarrow\infty,$%
\[
\frac{n}{2}\ln(2a)+\frac{1}{2}n\ln\left(  \frac{n}{a}\right)  +x\ln\left(
1-\sqrt{\frac{n}{a}}\right)  +\sqrt{an}-\sqrt{n}\sim\Phi_{3}\left(
\xi\right)  ,
\]
where%
\begin{equation}
\Phi_{3}\left(  \xi\right)  =\frac{n}{2}\ln\left(  2n\right)  -\frac{3}%
{2}n-\xi\sqrt{2n}. \label{phi3}%
\end{equation}
Also,%
\[
\sqrt{2\pi}\left(  \frac{n}{a}\right)  ^{\frac{1}{6}}\left(  \sqrt{a}-\sqrt
{n}\right)  ^{\frac{1}{3}}\sim n^{\frac{1}{6}}%
\]
and%
\[
\left(  \frac{n}{a}\right)  ^{\frac{1}{6}}\frac{\left(  \Omega^{-}-x\right)
}{\left(  \sqrt{a}-\sqrt{n}\right)  ^{\frac{2}{3}}}\sim-n^{\frac{1}{6}}%
\sqrt{2}\left(  \xi+\sqrt{2n}\right)  .
\]
Therefore,%
\begin{equation}
H_{n}(\xi)\sim\Lambda_{3}\left(  \xi\right)  \equiv(-1)^{n}\sqrt{2\pi}%
n^{\frac{1}{6}}\exp\left[  \Phi_{3}\left(  \xi\right)  \right]  \mathrm{Ai}%
\left[  -n^{\frac{1}{6}}\sqrt{2}\left(  \xi+\sqrt{2n}\right)  \right]  .
\label{H4}%
\end{equation}

\item Region V: $\xi\approx\sqrt{2n}.$

Using (\ref{eta}) in (\ref{c11}) we have, as $a\rightarrow\infty,$%
\[
\frac{n}{2}\ln(2a)+\frac{1}{2}n\ln\left(  \frac{n}{a}\right)  +x\ln\left(
1+\sqrt{\frac{n}{a}}\right)  -\sqrt{an}-\sqrt{n}\sim\Phi_{4}\left(
\xi\right)  ,
\]
with%
\begin{equation}
\Phi_{4}\left(  \xi\right)  =\frac{n}{2}\ln\left(  2n\right)  -\frac{3}%
{2}n-\xi\sqrt{2n}. \label{phi4}%
\end{equation}
Also,%
\[
\sqrt{2\pi}\left(  \frac{n}{a}\right)  ^{\frac{1}{6}}\left(  \sqrt{a}+\sqrt
{n}\right)  ^{\frac{1}{3}}\sim n^{\frac{1}{6}}%
\]
and%
\[
\left(  \frac{n}{a}\right)  ^{\frac{1}{6}}\frac{\left(  x-\Omega^{+}\right)
}{\left(  \sqrt{a}+\sqrt{n}\right)  ^{\frac{2}{3}}}\sim n^{\frac{1}{6}}%
\sqrt{2}\left(  \xi-\sqrt{2n}\right)  .
\]
Therefore,%
\begin{equation}
H_{n}(\xi)\sim\Lambda_{4}\left(  \xi\right)  \equiv\sqrt{2\pi}n^{\frac{1}{6}%
}\exp\left[  \Phi_{4}\left(  \xi\right)  \right]  \mathrm{Ai}\left[
n^{\frac{1}{6}}\sqrt{2}\left(  \xi-\sqrt{2n}\right)  \right]  . \label{h5}%
\end{equation}
Once again, we have
\[
\Lambda_{3}\left(  -\xi\right)  =(-1)^{n}\Lambda_{4}\left(  \xi\right)
\]

\item Region VI: $-\sqrt{2n}\ll\xi\ll\sqrt{2n}.$

From (\ref{c10}), we immediately obtain%
\[
H_{n}(\xi)\sim\Lambda_{5}\left(  \xi\right)  \equiv\Lambda_{1}\left(
\xi\right)  +\Lambda_{2}\left(  \xi\right)  .
\]
Since $-1<\frac{\xi}{\sqrt{2n}}<1,$ we set
\begin{equation}
\xi=\sqrt{2n}\sin\left(  \theta\right)  ,\quad-\frac{\pi}{2}<\theta<\frac{\pi
}{2}. \label{theta}%
\end{equation}
From (\ref{sigma}) we have%
\[
\sigma=\sqrt{2n}\cos\left(  \theta\right)  \mathrm{i.}%
\]
Thus,
\[
n\pi\mathrm{i}+\Phi_{1}\left(  \xi\right)  =\frac{n}{2}\left[  \ln\left(
2n\right)  -\cos\left(  2\theta\right)  \right]  +n\left[  \frac{1}{2}%
\sin\left(  2\theta\right)  +\theta-\frac{\pi}{2}\right]  \mathrm{i,}%
\]%
\[
\Phi_{2}\left(  \xi\right)  =\frac{n}{2}\left[  \ln\left(  2n\right)
-\cos\left(  2\theta\right)  \right]  -n\left[  \frac{1}{2}\sin\left(
2\theta\right)  +\theta-\frac{\pi}{2}\right]  \mathrm{i,}%
\]
and%
\[
U_{1}\left(  \xi\right)  =\frac{\exp\left(  \frac{\theta}{2}\mathrm{i}\right)
}{\sqrt{2\cos\left(  \theta\right)  }},\quad U_{2}\left(  \xi\right)
=\frac{\exp\left(  -\frac{\theta}{2}\mathrm{i}\right)  }{\sqrt{2\cos\left(
\theta\right)  }}.
\]
Hence,%
\begin{equation}
\Lambda_{5}\left[  \sqrt{2n}\sin\left(  \theta\right)  \right]  =\sqrt
{\frac{2}{\cos\left(  \theta\right)  }}\exp\left[  \Phi_{5}\left(
\theta\right)  \right]  \cos\left(  \Theta\right)  , \label{lambda5}%
\end{equation}
with%
\begin{equation}
\Phi_{5}\left(  \theta\right)  =\frac{n}{2}\left[  \ln\left(  2n\right)
-\cos\left(  2\theta\right)  \right]  \label{phi5}%
\end{equation}
and%
\begin{equation}
\Theta=n\left[  \frac{1}{2}\sin\left(  2\theta\right)  +\theta-\frac{\pi}%
{2}\right]  +\frac{\theta}{2}. \label{Theta}%
\end{equation}

Using (\ref{theta}), we can write (\ref{phi5}) and (\ref{Theta}) in terms of
$\xi$%
\[
\Phi_{5}\left(  \theta\right)  =\frac{n}{2}\left[  \ln\left(  2n\right)
-1\right]  +\frac{\xi^{2}}{2},
\]%
\[
\Theta=\frac{\xi}{2}\sqrt{2n-\xi^{2}}+\left(  n+\frac{1}{2}\right)
\arcsin\left(  \frac{\xi}{\sqrt{2n}}\right)  -n\frac{\pi}{2}.
\]
Also,%
\[
\sqrt{\frac{2}{\cos\left(  \theta\right)  }}=\sqrt{2}\left(  1-\frac{\xi^{2}%
}{2n}\right)  ^{-\frac{1}{4}}.
\]
Therefore,%
\begin{align}
\Lambda_{5}\left(  \xi\right)   &  =\sqrt{2}\left(  1-\frac{\xi^{2}}%
{2n}\right)  ^{-\frac{1}{4}}\exp\left\{  \frac{n}{2}\left[  \ln\left(
2n\right)  -1\right]  +\frac{\xi^{2}}{2}\right\} \label{lambda51}\\
&  \times\cos\left[  \frac{\xi}{2}\sqrt{2n-\xi^{2}}+\left(  n+\frac{1}%
{2}\right)  \arcsin\left(  \frac{\xi}{\sqrt{2n}}\right)  -n\frac{\pi}%
{2}\right]  .\nonumber
\end{align}
Considering the leading term of (\ref{lambda51}) as $n\rightarrow\infty,$ we
obtain%
\[
\Lambda_{5}\left(  \xi\right)  \sim\sqrt{2}\exp\left\{  \frac{n}{2}\left[
\ln\left(  2n\right)  -1\right]  +\frac{\xi^{2}}{2}\right\}  \cos\left(
n\frac{\pi}{2}-\xi\sqrt{2n}\right)
\]
in agreement with formula (4.14.9) in \cite{MR0350075}.
\end{enumerate}

\section{Zeros}

Let us denote by $\zeta_{1}^{n}>\zeta_{2}^{n}>\cdots>\zeta_{n}^{n}$ the zeros
of $H_{n}\left(  \xi\right)  ,$ enumerated in decreasing order. It then
follows from (\ref{lambda5}) that $\zeta_{j}^{n}=\sqrt{2n}\sin\left(  \tau
_{j}^{n}\right)  ,$ where $\tau_{j}^{n}$ satisfies%
\[
n\left[  \frac{1}{2}\sin\left(  2\tau_{j}^{n}\right)  +\tau_{j}^{n}-\frac{\pi
}{2}\right]  +\frac{\tau_{j}^{n}}{2}=\frac{\pi}{2}-j\pi,\quad1\leq j\leq n.
\]
We can rewrite the equation%
\[
n\left[  \frac{1}{2}\sin\left(  2t\right)  +t-\frac{\pi}{2}\right]  +\frac
{t}{2}=A
\]
as Kepler's equation%
\begin{equation}
E-\varepsilon\sin(E)=M,\label{Kepler}%
\end{equation}
with%
\begin{equation}
E=2t,\quad M=2\frac{2A+n\pi}{2n+1},\quad\varepsilon=-\frac{2n}{2n+1}.\label{M}%
\end{equation}
It is well known \cite{MR1349110} that the solution of (\ref{Kepler}) can be
expressed as a Kapteyn series%
\begin{equation}
E=M+2%
{\displaystyle\sum\limits_{k=1}^{\infty}}
\frac{1}{k}\mathrm{J}_{k}\left(  k\varepsilon\right)  \sin\left(  kM\right)
,\label{kapteyn}%
\end{equation}
where $\mathrm{J}_{k}\left(  \cdot\right)  $ is a Bessel function of the first kind.

Thus, using (\ref{M}) in (\ref{kapteyn}) with $A=\frac{\pi}{2}-j\pi,$ we
obtain%
\begin{equation}
\tau_{j}^{n}=\pi\frac{1+n-2j}{2n+1}+%
{\displaystyle\sum\limits_{k=1}^{\infty}}
\frac{1}{k}\mathrm{J}_{k}\left(  -\frac{2n}{2n+1}k\right)  \sin\left(
2\pi\frac{1+n-2j}{2n+1}k\right)  , \label{zeros}%
\end{equation}
for $1\leq j\leq n.$ Using the reflection formula \cite{atlas} $\mathrm{J}%
_{k}\left(  -x\right)  =\left(  -1\right)  ^{k}\mathrm{J}_{k}\left(  x\right)
,$ we can write (\ref{zeros}) as
\begin{equation}
\tau_{j}^{n}=\frac{\pi}{2}-\frac{\pi}{2}\left(  4j-1\right)  N^{-1}-%
{\displaystyle\sum\limits_{k=1}^{\infty}}
\frac{1}{k}\mathrm{J}_{k}\left[  \left(  1-N^{-1}\right)  k\right]
\sin\left(  \frac{4j-1}{N}k\pi\right)  , \label{zeros1}%
\end{equation}
where $N=2n+1.$




\end{document}